\documentclass{crm-article}
\newcommand{\argmin}{\mathop{\mathrm{argmin}}}
\def\dv#1{{\color{black}#1}} 

\begin{document}





\year=2021 

\journalSection{Mathematical modeling and numerical simulation}
\journalSectionEn{Mathematical modeling and numerical simulation}

\journalReceived{08.12.2020.}
\journalAccepted{20.12.2020.}

\UDC{519.8}
\title{Поиск равновесий в двухстадийных моделях распределения транспортных потоков по сети.}
\titleeng{Finding equilibrium in two-stage traffic assignment model.}

\author[1]{\firstname{Е.\,В.}~\surname{Котлярова}}
\authorfull{Котлярова Екатерина Владимировна}
\authoreng{\firstname{E.\,V.}~\surname{Kotliarova}}
\authorfulleng{Ekaterina V. Kotliarova}
\email{kotlyarova.ev@phystech.edu}

 \author[1, 2, 3]{\firstname{А.\,В.}~\surname{Гасников }}
 \authorfull{Александр Владимирович Гасников}
 \authoreng{\firstname{A.\,V.}~\surname{Gasnikov}}
 \authorfulleng{Alexander V. Gasnikov}
 \email{gasnikov@yandex.ru}
 
  \author[1]{\firstname{Е.\,В.}~\surname{Гасникова}}
 \authorfull{Евгения Владимировна Гасникова}
 \authoreng{\firstname{E.\,V.}~\surname{Gasnikova}}
 \authorfulleng{Evgenia V. Gasnikova}
 \email{egasnikov@yandex.ru}

  \author[1]{\firstname{Д.\,В.}~\surname{Ярмошик}}
 \authorfull{Демьян Валерьевич Ярмошик}
 \authoreng{\firstname{D.\,V.}~\surname{Yarmoshik}}
 \authorfulleng{Demyan V. Yarmoshik}
 \email{yarmoshik.dv@phystech.edu}

 \affiliation[1]{Национальный исследовательский университет «Московский физико-технический институт»,\protect\\ Россия, 141701, г. Долгопрудный, Институтский пер., д. 9}
 \affiliationeng{National Research University Moscow Institute of Physics and Technology,\protect\\ 9 Institute lane, Dolgoprudny, 141701, Russia}
  \affiliation[2]{Институт проблем передачи информации РАН,\protect\\ Россия, 127051, г. Москва, Б. Каретный пер., д. 9}
 \affiliationeng{Institute for Information Transmission Problems RAS,\protect\\ 9 B. Karetny lane, Moscow, 127051, Russia}
 \affiliation[3]{Кавказский математический центр, \protect\\ Россия, 385000, г. Майкоп, адрес Первомайская ул., 208}
 \affiliationeng{Caucasus Mathematical Center,\protect\\ 208 Pervomaiskaia street, Maikop, 385000, Russia}

\begin{abstract}
В работе описывается двухстадийная модель равновесного распределения транспортных потоков. Модель состоит из двух блоков. Первый блок -- модель расчёта матрицы корреспонденций, второй блок -- модель равновесного распределения транспортных потоков по путям. Первая модель, используя матрицу транспортных затрат (затраты на перемещение из одного района в другой район), рассчитывает матрицу корреспонденций, описывающих потребности в объемах передвижениях из одного района в в другой район. Вторая модель описывает на базе равновесного принципа Нэша--Вардропа (каждый водитель выбирает кратчайший для себя путь), как именно потребности в перемещениях, задаваемые матрицей корреспонденций, распределятся по возможным путям. Зная способы распределения потоков по путям, можно рассчитать матрицу затрат. Равновесием в двухстадийной модели транспортных потоков называют неподвижную точку цепочки из этих двух моделей. В статье предложен способ сведения задачи поиска описанного равновесия к задаче выпуклой негладкой оптимизации. Предложен численный способ решения полученной задачи оптимизации. Проведены численные эксперименты для небольших городов \dv{США и для Владивостока.} 
\end{abstract}

\keyword{модель расчета матрицы корреспонденций}
\keyword{многостадийная модель}
\keyword{модель равновесного распределения потоков по путям}

\begin{abstracteng}
Authors describe a two-stage traffic assignment model. It contains of two blocks. The first block consists of model for calculating  correspondence (demand) matrix, whereas the second block is a traffic assignment model. The first model calculates a matrix of correspondences using a matrix of transport costs. It characterizes the required volumes of movement from one area to another. The second model describes how exactly the needs for displacement, specified by the correspondence matrix, are distributed along the possible paths. It works on the basis of the Nash--Wardrop equilibrium (each driver chooses the shortest path). Knowing the ways of distribute flows along the paths, it is possible to calculate the cost matrix. Equilibrium in a two-stage model is a fixed point in the sequence of these two models. The article proposes a method of reducing the problem of finding the equilibrium to the problem of the convex non-smooth optimization. Also a numerical method for solving the obtained optimization problem is proposed. Numerical experiments were carried out for the small towns \dv{of USA and for Vladivostok.}
\end{abstracteng}

\keywordeng{correspondence matrix calculation model}
\keywordeng{multi stage model}
\keywordeng{equilibrium distribution model of traffic flow}

\maketitle

\paragraph{Введение}
В данной статье описывается (с обоснованием) вариационный (экстремальный) принцип, сводящий поиск равновесного распределения транспортных потоков по сети к задаче выпуклой оптимизации. Под распределением потоков понимается: 1) расчёт матрицы корреспонденций и 2) распределение потоков по путям при заданных корреспонденциях. Таким образом, речь идет о двухуровневой модели распределения. Многостадийные модели транспортных потоков являются одним из основных объектов изучения при долгосрочном транспортном планировании \cite{Ortuzar2002,gasnikov2013book,gasnikov2020posobie}. С помощью таких моделей можно просчитывать долгосрочные последствия: введения в эксплутатацию различных инфраструктурных объектов, изменения дорожной сети и т.п.

Следуя работам \cite{gasnikov2014matmod,babicheva2015mipt,gasnikov2016diss} в статье выписывается задача выпуклой оптимизации, к которой сводится поиск равновесия в такой двухстадийной модели. Далее эта задача упрощается (путем перехода к двойственному представлению), и описывается численный способ решения возникающей в итоге (двойственной) задачи. Отличительными особенностями данной работы являются: 1) простой способ получения итоговой задачи оптимизации и способа ее решения; 2) проведенные численные эксперименты. За базу был взят код  \cite{kubentaeva2020code}, в котором  рассматривалось только распределение потоков по путям при заданных корреспонденциях \cite{gasnikov2018kim,baimurzina2019jvm,Kubentayeva2020}. В качестве источника данных использовался ресурс \cite{Stabler2020}. \dv{Также из разных источников были собраны данные по Владивостоку.} Таким образом в данной статье код \cite{kubentaeva2020code} был распространен на поиск равновесий в двухстадийных моделях, введенных в \cite{gasnikov2014matmod,babicheva2015mipt,gasnikov2016diss}.

\paragraph{Основные определения и обозначения}
Для простоты будем рассматривать замкнутую транспортную систему, описываемую графом $G = \langle V, E\rangle$, где $V$ -- множество вершин ($|V| = n$), а $E$ -- множество ребер ($|E| = m$). Будем обозначать ребра графа через $e\in E$. Для стандартной транспортной системы можно ожидать, что $m\simeq 3n$. Для больших мегаполисов (таких как Москва)  $n \simeq 10^5$. 
Однако в данной работе мы будем рассматривать в основном примеры $n \lesssim 10^4$. Транспортный граф $G$ считается известным.

Часть вершин $O\subseteq V$ (\textit{origin}) являются источниками корреспонденций, а часть стоками корреспонденций $D\subseteq V$ (\textit{destination}). Если говорить более точно, то вводится множество пар (источник, сток) корреспонденций $OD \subseteq V\bigotimes V$. Сами корреспонденции будем обозначать через $d_{ij}$, где $(i,j)\in OD$. Как правило $|OD|\ll n^2$ \cite{gasnikov2014matmod}. Не ограничивая общности, будем далее считать, что $\sum_{(i,j)\in OD} d_{ij} = 1$. Множество пар $OD$ считается известным. Корреспонденции -- не известны! Однако известны (заданы) характеристики источников и стоков корреспонденций. То есть известны величины $\{l_i\}_{i\in O}$, $\{w_j\}_{j\in D}$
\begin{equation}\label{corr}
    \sum_{j: (i,j)\in OD} d_{ij} = l_i, \quad \sum_{i: (i,j)\in OD} d_{ij} = w_j.
\end{equation}
Заметим, что $\sum_{i\in O} l_i = \sum_{j\in D} w_j = 1$. Условие \eqref{corr} будем также для краткости записывать в виде $d\in (l,w)$.

Обозначим через $\tau_e(f_e)$ -- функцию затрат (например, временных) на проезд по ребру (участку дороги) $e$, если поток автомобилей на этом участке $f_e$. Функции $\tau_e(f_e)$ считаются заданными, например, таким образом: \cite{gasnikov2013book,Patriksson2015,gasnikov2020posobie}
\begin{equation}\label{BPR}
    \tau_e(f_e) = \bar{t}_e\left(1 +\kappa\left(\frac{f_e}{\bar{f}_e}\right)\right)^{\frac{1}{\mu}}, 
\end{equation}
где $\bar{t}_e$ -- время прохождения ребра $e$, когда участок свободный (определяется разрешенной скоростью на данном участке), а $\bar{f}_e$ -- пропускная способность ребра $e$ (определяется полосностью: [пропускная способность] $\le$ [число полос] * [2000 авт/час] и характерстиками перекрестков). Считается, что эти характеристики известны \cite{Stabler2020}. Параметр $\mu = 0.25$ -- BPR-функции \cite{Patriksson2015}, но допускается и $\mu\to 0+$ -- модель стабильной динамики \cite{Nesterov2003,gasnikov2013book,gasnikov2014matmod,gasnikov2016matmod,gasnikov2016diss,Gasnikov2018,gasnikov2020posobie}. Параметр $\kappa >0$ также считается заданным.

Полезно также ввести $t_e$ -- (временные) затраты на прохождения ребра $e$. Согласно вышенаписанному $t_e = \tau_e(f_e)$. По этим затратам $t = \{t_e\}_{e\in E}$ можно определить затраты на перемещение из источника $i$ в сток $j$ по кратчайшему пути: $T_{ij}(t) = \min_{p \in P_{ij}} \sum_{e\in E} \delta_{ep}t_e$, где $p$ -- путь (без самопересечений -- циклов) на графе (набор ребер), $P_{ij}$ -- множество всевозможных путей на графе, стартующих из источника $i$ и заканчивающихся в стоке $j$, $\delta_{ep} = 1$, если ребро $e$ принадлежит пути $p$ и $\delta_{ep} = 1$ -- иначе.

В ряде выкладок далее также будет полезен вектор $x = \{x_p\}_{p \in P}$ -- вектор распределения потоков по путям, где $P = \bigcup_{(i,j)\in OD} P_{ij}$. Заметим, что $f_e = \sum_{p} \delta_{ep} x_p$ или в матричном виде $f = \Theta x$, где $\Theta = \|\delta_{ep}\|_{e\in E, p\in P}$.

\paragraph{Энтропийная модель расчёта матрицы корреспонденицй}

Под энтропийной моделью расчета матрицы корреспонденций $d(T)$ понимается определенный способ вычисления набора корреспонденций $\{d_{ij}\}_{(i,j)\in OD}$ по известной матрице затрат $\{T_{ij}\}_{(i,j)\in OD}$. Этот способ заключается в решении задачи энтропийно-линейного программирования, которую можно понимать, как энтропийно-регуляризованную транспортную задачу\footnote{Вместо $d_{ij}\ln d_{ij}$ точнее было бы писать $d_{ij}\ln\left(d_{ij}/\left(\sum_{(i,j)\in OD}d_{ij}\right)\right)$ \cite{gasnikov2014matmod}, но $\sum_{(i,j)\in OD} d_{ij} = 1$, поэтому, возможно, упрощенная форма записи.}
\begin{equation}\label{Wilson}
   \min_{d\in (l,w);d\ge 0} \sum_{(i,j)\in OD} d_{ij}T_{ij} +  \gamma \sum_{(i,j)\in OD} d_{ij}\ln d_{ij},
\end{equation}
где параметер $\gamma >0$ считается известным \cite{Wilson1978,gasnikov2010,gasnikov2013book,gasnikov2016,gasnikov2016diss,gasnikov2020posobie}. Относительно выбора этого парметра, см. \cite{gasnikov2014matmod,gasnikov2020posobie,ivanova2020}.

\paragraph{Модели равновесного распределения транспортных потоков по путям}
Матрица корреспонденций $\{d_{ij}\}_{(i,j)\in OD}$ порождает (вообще говоря, неоднозначно) некий вектор распределения потоков по путям $x$. Неоднозначность заключается в том, что балансовые ограничения, которые возникают на $x \in X(d)$:
\begin{equation*}\label{Xd}
  x\ge 0:\quad \forall (i,j)\in OD \to  \sum_{p\in P_{ij}} x_p = d_{ij},
\end{equation*}
как правило, не определяют вектор $x$ однозначно. Вектор $x$, в свою очередь, пораждает вектор потоков на ребрах, $f = \Theta x$, который, в свою очередь, порождает вектор (временных) затрат на ребрах $t(f) = \{\tau_e(f_e)\}_{e\in E}$. На основе последнего вектора уже можно рассчитать матрицу затрат на кратчайших путях $T(t) = \{T_{ij}(t)\}_{(i,j)\in OD}$. Собственно, модель равновесного распределения потоков это формализация \textit{принципа Нэша--Вардропа} о том, что в равновесии каждый водитель выбирает для себя кратчайший путь \cite{gasnikov2013book,Patriksson2015,gasnikov2020posobie}. Другими словами, если для заданной корреспонденции $(i,j)\in OD$ известно, что  (\textit{условие комплиментарности})
\begin{center}
   $x_{p'} >0$, где $p' \in P_{ij}$, то $T_{ij}(t) = \min_{p \in P_{ij}} \sum_{e\in E} \delta_{ep}t_e = \sum_{e\in E} \delta_{ep'}t_e$. 
\end{center}
Задача поиска равновесия сводится, таким образом, к поиску такого вектора $x \in X(d)$, который бы порождал такие затраты $T := T(t(f(x)))$, что выполянется условие комплиментарности. В написанном выше виде искать равновесный вектор $x \in X(d)$ представляется сложной задачей, сводящейся к решению системы нелинейных уравнений. Однако, в данном случае (рассматривается потенциальная игра загрузки) можно свести поиск равновесия к решению задачи выпуклой оптимизации\footnote{Подобно \eqref{Wilson} можно искать не равновесия Нэша--Вардропа, а стохастическое равновесия. Это приводит к дополнительному энтропийному слагаемому в \eqref{Beckman} \cite{gasnikov2013book,gasnikov2014matmod,gasnikov2015MIPTershov,baimurzina2019jvm,gasnikov2020posobie}.}
\begin{equation}\label{Beckman}
   \min_{(f,x): f=\Theta x; x\in X(d)} \sum_{e\in E} \int_{0}^{f_e} \tau_e(z) dz.
\end{equation}
Решение задачи дает модель вычисления вектора потока на ребрах при заданной матрице корреспонденций $f(d)$  \cite{gasnikov2013book,gasnikov2014matmod,Patriksson2015,gasnikov2016matmod,gasnikov2020posobie}.  

\paragraph{Двухстадийная модель}
Выше были описаны две модели. В первой (расчёт матрицы корреспонденций) на вход подается матрица затрат $T$, а на выходе получается матрица корреспонденций  $d(T)$. Во второй модели наоборот, на вход подается матрица корреспонденций  $d(T)$, а на выходе рассчитывается матрица затрат $T(d) = T(t(f(d)))$. 

Под равновесием в двухстадийной транспортной модели понимается такая пара $(f,d)$, что \cite{Ortuzar2002,gasnikov2014matmod,babicheva2015mipt,gasnikov2016diss,gasnikov2020posobie}
\begin{equation}\label{FP}
d = d(T(t(f))), \quad f = f(d),
\end{equation}
то есть $(f,d)$ -- есть неподвижная точка описанных двух блоков моделей. Собственно, часто на практике так и ищут равновесие последовательно (друг за другом) прогоняя описанные два блока \cite{Ortuzar2002,gasnikov2014matmod}. Однако, насколько нам известно, нет никаких теоретических гарантий, что такая процедура (последовательная прогонка) будет сходиться к неподвижной точке. Собственно, описанные в следующих разделах численные эксперименты показывают, что на практике сходимость наблюдается далеко не всегда. Но даже если наблюдается сходимость, то непонятно, насколько эта сходимость может быть быстрой, и лучший ли это способ (простая прогонка) численного решения \eqref{FP}? Далее, следуя \cite{gasnikov2014matmod,gasnikov2015usik,babicheva2015mipt,gasnikov2016diss,gasnikov2020posobie}, предлагается эквивалентный способ перезаписи задачи \eqref{FP} как задачи выпуклой оптимизации, которую можно уже  решать оптимальными по скорости (глобально сходящимися) алгоритмам.

\begin{teo}
Задача поиска неподвижной точки \eqref{FP} сводится к задаче выпуклой оптимизации
\begin{equation}\label{TS}
\min_{(f,x,d): f=\Theta x; x\in X(d); d\in(l,w); d\ge 0} \sum_{e\in E} \int_{0}^{f_e} \tau_e(z) dz +  \gamma \sum_{(i,j)\in OD} d_{ij}\ln d_{ij} .
\end{equation}
\end{teo} 

\textbf{Схема доказательства}.
Вернёмся к формулировке задачи оптимизации \eqref{Beckman} и заметим следующее важное свойство функционала задачи $\Psi(x) = \sum_{e\in E} \int_{0}^{f_e(x)} \tau_e(z) dz$: для всех $p \in P$
$$\frac{\partial \Psi(x)}{\partial x_p} = T_p(x),$$
где $T_p(x) = \sum_{e\in E}\delta_{ep}\tau_e(f_e(x))$ -- затраты на пути $p$. Это и есть проявление того, что рассматривается потенциальная игра загрузки. Собственно, для  \eqref{Wilson} также можно выделить потенциал $\Phi(d) = \sum_{(i,j)\in OD} d_{ij}T_{ij}$, для которого: для всех $(i,j)\in OD$
$$\frac{\partial \Phi(d)}{\partial d_{ij}} = T_{ij}.$$
Если бы удалось найти такую функцию $\tilde{\Phi}(d)$, для  которой: для всех $(i,j)\in OD$
\begin{equation}\label{PD}
    \frac{\partial \tilde{\Phi}(d)}{\partial d_{ij}} = T_{ij}(t(f(d))),
\end{equation}
где $T_{ij}(t(f(d)))$ определяется из решения задачи \eqref{Beckman}, то решение задачи 
 $$\min_{d\in (l,w);d\ge 0} \tilde{\Phi}(d) +  \gamma \sum_{(i,j)\in OD} d_{ij}\ln d_{ij}$$
 давало бы равновесную матрицу корреспонденций $d$, по которой можно было бы уже оценить и равновесный вектор потоков на ребрах $f = f(p)$, см. \eqref{FP}. Детали см.  \cite{gasnikov2014matmod,gasnikov2015usik,babicheva2015mipt,gasnikov2016diss,gasnikov2020posobie}. 
 
 Таким образом, цель -- найти потецниал $\tilde{\Phi}(d)$, если, конечно,  он существует. Покажем, что $$\tilde{\Phi}(d) = \min_{(f,x): f=\Theta x; x\in X(d)} \sum_{e\in E} \int_{0}^{f_e} \tau_e(z) dz,$$
 см. формулу \eqref{Beckman}.
 
 Для этого введем (выпуклые) функции $\sigma_e(f_e) = \int_{0}^{f_e} \tau_e(z) dz$, и обозначим сопряженные к ним функции через $\sigma_e^*(t_e) = \max_{f_e \ge 0} \left\{f_e t_e - \sigma_e(f_e) \right\}$. Тогда (детали см. в \cite{gasnikov2014matmod,gasnikov2016diss,gasnikov2020posobie})
 $$\tilde{\Phi}(d) =\min_{(f,x): f=\Theta x; x\in X(d)} \sum_{e\in E} \sigma_e(f_e) =\min_{(f,x): f=\Theta x; x\in X(d)} \sum_{e\in E} \max_{t_e\in \text{dom }\sigma_e^*}\{f_e t_e - \sigma_e^*(t_e)\} = $$
 $$\max_{t_e\in \text{dom }\sigma_e^*, e \in E} \left\{\min_{(f,x): f=\Theta x; x\in X(d)} \sum_{e\in E} f_e t_e\right\} - \sum_{e\in E} \sigma_e^*(t_e) =$$
 \begin{equation}\label{dual}
    \max_{t_e\in \text{dom }\sigma_e^*, e \in E} \sum_{(i,j)\in OD} d_{ij} T_{ij}(t) - \sum_{e\in E}  \sigma_e^*(t_e). 
 \end{equation}
 Здесь $\text{dom }\sigma_e^*$ означает область определения функции $\sigma_e^*(t_e)$. Из формулы \eqref{dual} и формулы Демьянова--Данскина \cite{Bertsekas2009,gasnikov2021book} следует \eqref{PD}.
 
 Строго говоря, приведенные выше рассуждения еще не являются доказательством, поскольку апеллируют к понятию потенциала и использованию соответствующших теорем популяционной теории игр загрузки. Однако вместо того, чтобы приводить здесь эти соображения (подобно тому, как это было сделано, например, в работах \cite{Sandholm2010,gasnikov2013book,gasnikov2014matmod,gasnikov2015usik,gasnikov2016,Dvurechensky2016,gasnikov2020posobie}), в этой статье мы ограничимся тем, что заметим, что выписанная задача \eqref{TS}, как задача оптимизации относительно $d$ при <<заморожоженных>> $(f,x)$, совпадает с задачей \eqref{Wilson} и, наоборот, при <<замороженном>> $d$ задача \eqref{TS}, как задача оптимизации относительно $(f,x)$, совпадает с задачей \eqref{Beckman}. Таким образом, если удалось найти такую задачу, решение которой одновременно дает нужные нам связи переменных, описываемые формулой \eqref{FP}, то это и означает, что нам удалось свести поиск неподвижной точки сложного нелинейного отображения (которое не удается выписать аналитически) к явно выписанной задаче оптимизации \eqref{TS}. Ометим, что решение  этой выпуклой задачи оптимизации  по сложности сопоставимо с решением задачи \eqref{Beckman}, что будет пояснено в следующем разделе.

\paragraph{Переход к двойственной задаче}
Как следует из схемы доказательства теоремы 1, задача выпуклой оптимизации \eqref{TS} можно переписать эквивалентным седловым образом, введя двойственные переменные $t = \{t_e\}_{e\in E}$,\footnote{Отметим, что, как уже подмечалось ранее, есть явная связь двойственных переменных с прямыми:  $t_e = \tau_e(f_e)$, $e\in E$. При $\mu\to 0+$ (см. \eqref{BPR}) эта связь становится более хитрой, см. \cite{gasnikov2014matmod,gasnikov2016matmod,Gasnikov2018,baimurzina2019jvm,Kubentayeva2020,gasnikov2020posobie}. Это же замечание имеет место и для формулы \eqref{grad} далее.} которые имеют естественную интерпретацию вектора потоков на ребрах, 
$$ \min_{d\in (l,w);d\ge 0}  \max_{t_e\in \text{dom }\sigma_e^*, e \in E} \left\{ \sum_{(i,j)\in OD} d_{ij} T_{ij}(t) - \sum_{e\in E} \sigma_e^*(t_e) + \gamma \sum_{(i,j)\in OD} d_{ij}\ln d_{ij}\right\}.$$
Последнюю задачу удобнее переписать в виде\footnote{Это седловая задача. Отметим, что известные сейчас оптимальные методы решения выпукло-вогнутых седловых задач (см., например, \cite{gasnikov2021book,gasnikov2021} и цитированную там литературу) здесь не подходят, поскольку не получается эффективно проектироваться на ограничение $d\in(l,w)$. Поэтому далее это ограничение заносится в функционал с помощью принципа множителей Лагранжа.}
$$
  \max_{t_e\in \text{dom }\sigma_e^*, e \in E}\quad \min_{d\in (l,w);\sum_{(i,j)\in OD}d_{ij}=1;d\ge 0}  \left\{ \sum_{(i,j)\in OD} d_{ij} T_{ij}(t) + \gamma \sum_{(i,j)\in OD} d_{ij}\ln d_{ij}\right\} - \sum_{e\in E} \sigma_e^*(t_e).
$$
Вспомогательную задачу минимизации можно представить через двойственную к ней:
\begin{equation*}
 \max_{t_e\in \text{dom }\sigma_e^*, e \in E; (\lambda, \mu)} -\gamma\ln\left(\sum_{(i,j)\in OD}\exp\left(\frac{-T_{ij}(t) + \lambda_i + \mu_j}{\gamma}\right)\right) + \langle l,\lambda \rangle + \langle w,\mu \rangle - \sum_{e\in E} \sigma_e^*(t_e)= 
\end{equation*}
\begin{equation}\label{TSD}
 -\min_{t_e\in \text{dom }\sigma_e^*, e \in E; (\lambda, \mu)} \gamma\ln\left(\sum_{(i,j)\in OD}\exp\left(\frac{-T_{ij}(t) + \lambda_i + \mu_j}{\gamma}\right)\right) - \langle l,\lambda \rangle - \langle w,\mu \rangle + \sum_{e\in E} \sigma_e^*(t_e). 
\end{equation}
Обратим внимание, что добавленное по $d$ ограничение $\sum_{(i,j)\in OD}d_{ij}=1$ тавтологично, поскольку следует из $d\in(l,w)$. Тем не менее, удобнее его добавить, чтобы при взятии $\min$ получалась равномерно гладкая функция (типа softmax), а не сумма экспонент, имеющая неограниченные константы гладкости \cite{gasnikov2015MIPTnesterov,gasnikov2020posobie}. Множители $\lambda$ и $\mu$ являются двойственными множителями (множителями Лагранжа) к ограничениям $d\in(l,w)$ (см. \eqref{corr}), которые заносятся в функционал (ограничения $\sum_{(i,j)\in OD}d_{ij}=1;d\ge 0$ не заносятся в функционал). Заметим, что если $(t,\lambda,\mu)$ -- решение задачи \eqref{TSD}, то\footnote{Здесь $C_{\lambda}$ и $C_{\mu}$ -- произвольные числа.} $\left(t,\lambda + (C_{\lambda},...,C_{\lambda})^T,\mu+(C_{\mu},...,C_{\mu})^T\right)$ -- также будет решением задачи, т.е. решение задачи \eqref{TSD} не единственное \cite{gasnikov2015MIPTnesterov}. Заметим также, что зная $(\lambda,\mu)$, можно посчитать матрицу корреспонденций \cite{gasnikov2015MIPTnesterov,gasnikov2016jvm,Dvurechensky2020,gasnikov2020posobie}:
\begin{equation}\label{d}
    d_{ij}(\lambda,\mu)=\frac{\exp\left(\frac{-T_{ij}(t) + \lambda_i + \mu_j}{\gamma}\right)}{\sum_{(k,l)\in OD}\exp\left(\frac{-T_{kl}(t) + \lambda_k + \mu_l}{\gamma}\right)}.
\end{equation}
Для решения задачи выпуклой оптимизации (но, вообще говоря, не гладкой, поскольку функции $T_{ij}(t)$ -- негладкие)  можно использовать субградиентные методы. А именно, субградиент (далее обозначаем (супер-)субградиент таким же символом, как и градиент $\nabla$) целевого функционала по $t$ (стоящего под минимумом) \eqref{TSD} можно посчитать по формуле Демьянова--Данскина \cite{Bertsekas2009,gasnikov2021book}:
\begin{equation}\label{grad}
   \sum_{(i,j)\in OD} d_{ij}(\lambda,\mu)\nabla T_{ij}(t) + f = \sum_{(i,j)\in OD} d_{ij}(\lambda,\mu)\nabla T_{ij}(t) +  \left(\{\tau_e^{-1}(t_e)\}_{e\in E}\right)^T,
\end{equation}
где $\tau_e^{-1}$ -- обратная функция к $\tau_e$.
Примечательно, что отличие формулы \eqref{grad} от ее аналога, который можно получить, решая задачу \eqref{Beckman} посредством перехода к двойственной задаче  \cite{gasnikov2016matmod,gasnikov2016diss,Gasnikov2018,gasnikov2020posobie,Kubentayeva2020} только в том, что $d_{ij} = d_{ij}(\lambda,\mu)$, где $(\lambda,\mu)$ определяются из решения задачи:
\begin{equation}\label{EntrDual}
 \min_{(\lambda, \mu)} D(t,\lambda,\mu):= \gamma\ln\left(\sum_{(i,j)\in OD}\exp\left(\frac{-T_{ij}(t) + \lambda_i + \mu_j}{\gamma}\right)\right) - \langle l,\lambda \rangle - \langle w,\mu \rangle.
\end{equation}
Важное наблюдение, сделанное в работах \cite{gasnikov2015matmod,babicheva2015mipt,gasnikov2015MIPTnesterov,gasnikov2016diss,gasnikov2020posobie} заключается в том, что решать задачу \eqref{TSD} выгоднее не как задачу оптимизации по переменным $(t,\lambda,\mu)$, а как задачу только по переменной $t$, в то время как переменные $(\lambda,\mu)$ лишь используются для подсчета субградиента целевого функционала по формуле \eqref{grad} с  $d_{ij}(\lambda,\mu)$ рассчитываемым по формуле \eqref{d}, в которой

\begin{equation}\label{lm}
    (\lambda, \mu):= (\lambda(t), \mu(t)) \in \argmin_{(\lambda, \mu)} D(t,\lambda,\mu)  
\end{equation}

определяются как решение задачи \eqref{EntrDual}. Заметим, что сложность вычисления $\nabla T_{ij} (t)$ оптимальным алгоритмом Дейкстры (детали см., например, \cite{gasnikov2016matmod,Gasnikov2018,gasnikov2020posobie}) будет сопоставима со сложностью вычисления (с нужной  точностью) матрицы $d_{ij}(\lambda(t),\mu(t))$ \cite{Dvurechensky2018,Peyre2019,Guminov2019,Stonyakin2020,Tupitsa2020,gasnikov2021book}. Таким образом, получается, что сложность вычисления субградиента для двойственной задачи к \eqref{Beckman} и для задачи \eqref{TSD} сопоставимы. При этом, свойства гладкости (определяющие скорость сходимости используемых методов) целевого функционала в задаче \eqref{TSD} при переходе от оптимизации в пространстве $(t,\lambda,\mu)$  к оптимизации по переменной $t$ могут только улучшиться \cite{gasnikov2015MIPTnesterov,gasnikov2020posobie} (во всяком случае не ухудшиться). 

\begin{figure}[!ht]
    \centering
    \includegraphics[width=0.9\textwidth]{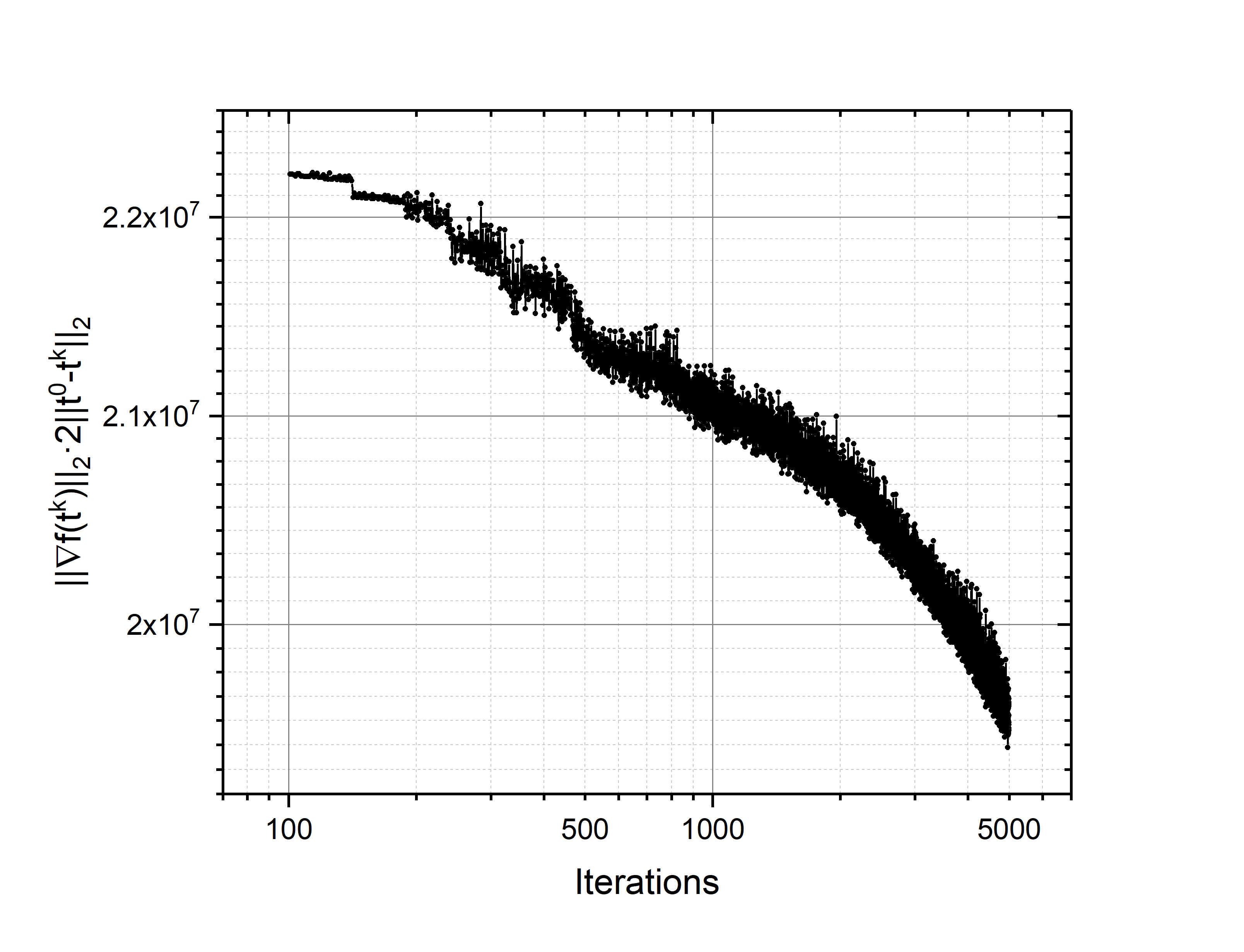}
    \caption{По оси абсцисс итерации, по оси ординат величина, оценивающая <<зазор двойственности>> $\Delta(t^k,B_k) = \max_{t \in B_k \cap \text{ dom }\sigma^*} \langle \nabla F(t^k), t^k  - t \rangle$}
    \label{fig:result}
\end{figure}

\paragraph{Вычислительные эксперименты}
 Итак, в качестве задачи оптимизации предлагается решать задачу: 
 \begin{equation}\label{f}
   \min_{t_e\in \text{dom }\sigma_e^*, e \in E} F(t):=  D(t,\lambda(t),\mu(t)) + \sum_{e\in E} \sigma_e^*(t_e).
 \end{equation}
 При этом, согласно \eqref{grad},
 $$\nabla F(t) = \sum_{(i,j)\in OD} d_{ij}(\lambda(t),\mu(t))\nabla T_{ij}(t) + \left(\{\tau_e^{-1}(t_e)\}_{e\in E}\right)^T,$$
 где $d_{ij}(\lambda(t),\mu(t))$ определяется формулами \eqref{d}, \eqref{EntrDual}, \eqref{lm}.
 В качестве способа решения задачи \eqref{f} предлагается использовать универсальный ускоренный градиентный метод, адаптивно настраивающийся на гладкость задачи \cite{Nesterov2015,gasnikov2015MIPTnesterov,gasnikov2015usik,Gasnikov2018,baimurzina2019jvm,Kamzolov2020,Nesterov2020,gasnikov2020posobie}. В худшем случае можно ожидать такой $F(t^k) - F(t_*) \simeq  O(k^{-1/2})$ скорости сходимости, где $t_*$ -- решение задачи \eqref{f} (эта скорость сходимости в общем случае не улучшаема для данного класса задач, см., например, \cite{gasnikov2021book} и цитированную там литературу), однако, как следует из результатов \cite{baimurzina2019jvm} на практике можно ожидать $F(t^k) - F(t_*) \simeq  O(k^{-1})$. Поскольку оптимальное значение целевого функционала $F(t_*)$ типично недоступно, то в качестве критерия оценки скорости сходимости  в экспериментах планиурется использовать величину, оценивающую <<зазор двойственности>> $\Delta(t^k,B_k) = \max_{t \in B_k \cap \text{ dom }\sigma^*} \langle \nabla F(t^k), t^k  - t \rangle$ (здесь для простоты считаем, что $\text{ dom }\sigma^*_e \equiv \text{ dom }\sigma^*$). Заметим, что с одной стороны $F(t^k) - F(t_*) \le \Delta(t^k,B_k)$, если $t_* \in B_k$, а с другой стороны, оценки скорости сходимости (в том числе отмеченные выше) получены, в действительности, для зазора двойственности, поскольку универсальный ускоренный градиентный метод является прямо-двойственным методом \cite{Nesterov2015,gasnikov2015MIPTnesterov,gasnikov2015usik,Gasnikov2018,baimurzina2019jvm,Kamzolov2020,Nesterov2020,gasnikov2020posobie,gasnikov2021book}. В численных экспериментах, базируясь на результатах о локализации последовательности, генерируемой методами типа (ускоренного) градиентного спуска  \cite{Gasnikov2018,gasnikov2021book}, в качестве множества $B_k$ выбирается евклидов шар с центром в точке $t_k$ и радиусом равным $2\|t^0 - t^k\|_2$.
 Для решения вспомогательной задачи \eqref{lm} можно использовать алгоритм Синхорна \cite{Dvurechensky2018,Peyre2019}, который в транспортной литературе чаще называют методом балансировки или методом Брэгмана(--Шелейховского)\cite{Wilson1978,gasnikov2013book}. Можно использовать и ускоренные варианты этого метода \cite{Guminov2019,Tupitsa2020}. Важной особенностью этой линейки методов (альтернативных направлений) является высокая (линейная) скорость сходимости при немалых значениях $\gamma > 0$. К сожалению, точной теории, насколько нам известно, этого режима сходимости еще нет, однако есть подтверждающие отмеченное наблюдение численные эксперименты и некоторые объяснения почему такая сходимость может быть \cite{Kroshnin2019,Stonyakin2020}. Стоит отметить, что использованние данного подхода предполагает, что для всех $i\in O$, $j \in D$ выполянется: $l_i,w_j >0$ и $T_{ij} < \infty$. В общем случае $OD \neq O\otimes D$ и какие-то $T_{ij} = \infty$ (что означает невозможность добраться из района $i$ в район $j$), также как, возможны районы, которые являются только жилыми (рабочими) районами, что приводит к $l_i = 0$ ($w_j = 0$).  Последнюю проблему можно решать чисто техническим образом, см., например, \cite{Dvurechensky2018,Kroshnin2019}, перераспределяя немного (в зависимости от желаемой точности решения задачи) с больших компонент этих векторов на нулевые, сохраняя нормировку. Проблему $T_{ij} = \infty$ также можно решать искусственно вводя пути из $i$ в $j$ приписывая им большие (но конечные) затраты (все это можно и содержательно проинтерпретировать). Для решения задачи \eqref{f} использовался вариант ускоренного универсального градиентного метода Нестерова \cite{Nesterov2015}, построенный на базе ускоренного (быстрого) градиентного метода подобных треугольников \cite{gasnikov2018jvm}. Код метода (и его адаптация к задаче \eqref{f}) был взят из \cite{baimurzina2019jvm,kubentaeva2020code,Kubentayeva2020}. В качестве точки старта метода выбиралось значение $t^0 = \bar{t}$, см. \eqref{BPR}.  В качестве источника данных ($G,O,D,OD,l,w,\Theta,\bar{t},\bar{f},\gamma,\kappa,\mu$) выбирались различные небольшие города (например, 
Су Фолс \cite{SiouxFalls}) с ресурса \cite{Stabler2020}. Проводилось два типа экспериментов. В первом подходе использовался метод простой (прямой) прогонки двух блоков (расчёт матрицы корреспонденций алгоритмом Синхорна, затем вычисление равновесного распределения потоков по путям универсальным ускоренным методом подобных треугольников в реализации \cite{baimurzina2019jvm,kubentaeva2020code,Kubentayeva2020}, потом пересчет матрицы затрат и снова вычисление матрицы корреспонденций и т.д.). Второй подход базировался на сочетании универсального ускоренного метода подобных треугольников в реализации \cite{baimurzina2019jvm,kubentaeva2020code,Kubentayeva2020} с алгоримтом Синхорна, для расчёта матрицы корреспонденций так как было описано выше в этой статье. На небольшом городе Су Фолс \cite{SiouxFalls} по модели Бэкмана с $\mu = 0.25$ не удалось добиться сходимости первого подхода. Результаты вычислительного эксперимента по второй модели приведены на рис.~\ref{fig:result}.

Из графика видно, что 
$$F(t^k) - F(t_*) \le \Delta(t^k,B_k) \simeq  O(k^{-1.67}),$$
что почти соответствует скорости сходимости в гладком случае $\sim k^{-2}$ \cite{gasnikov2021book}. Подчеркнем, что рассматриваемая задача \eqref{f} существенно негладкая, поэтому наблюдаемый результат можно интерпретировать таким образом, что негладкость задачи при правильном взгляде на нее (сквозь призму универсальных методов, настраивающихся на гладкость) не играет существенной роли в сложности ее численного решения.

Более подробно результаты экспериментов можно посмотреть по ссылке \cite{kоtlyarova2020code}.

\begin{figure}[!ht]
    \centering
    \includegraphics[width=0.9\textwidth]{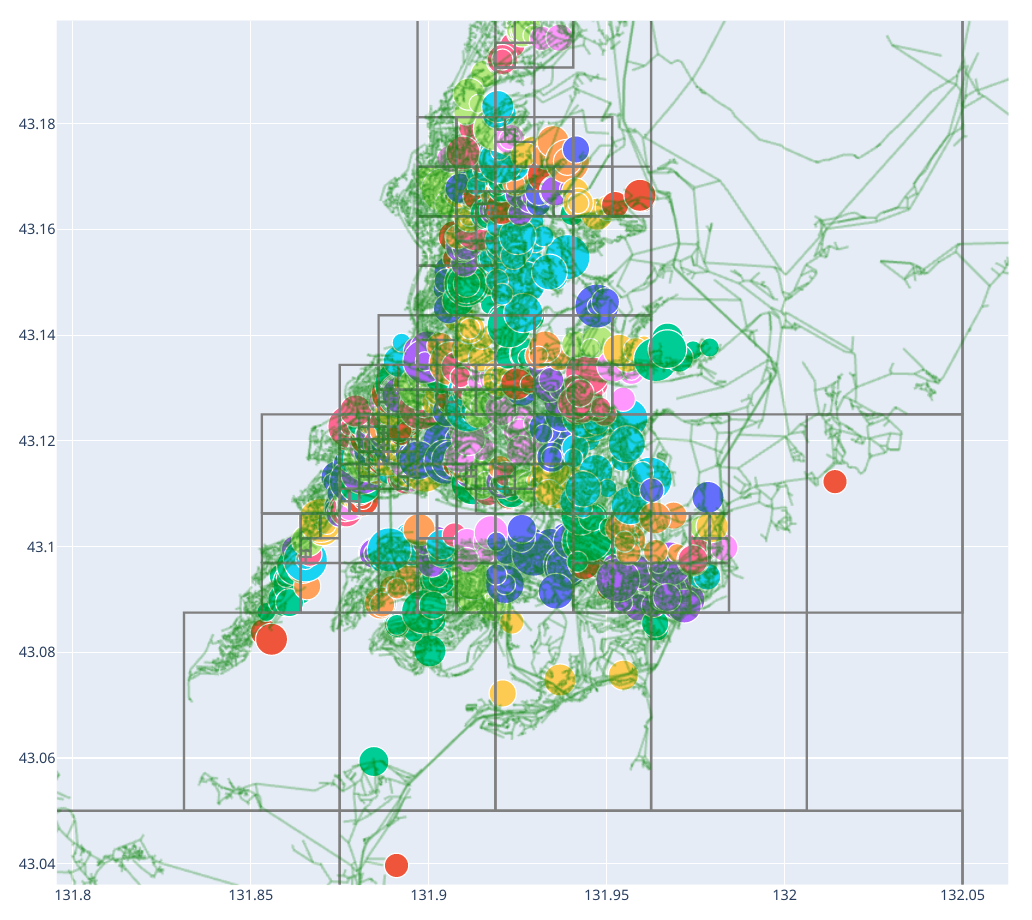}
    \caption{\dv{Транспортная сеть Владивостока. По осям --- географические координаты. Серые линии --- сетка разбиения на районы. Кружки --- предприятия, радиус соответствует логарифму числа рабочих мест, цвет --- принадлежности району.}}
    \label{fig:network}
\end{figure}

\dv{
Также был произведён расчёт двухстадийной модели для транспортной сети Владивостока.
Дорожный граф и данные о расположении мест жительства были предоставлены Е.А. Нурминским, а данные об адресах и количествах рабочих мест по ним были собраны с ресурса \cite{MinTrudPrim2020data}. 
Адреса для последующей обработки были преобразованы в географические координаты с помощью сервиса \cite{Dadata2020data}.
Для выбора источников и стоков корреспонденций и определения их характеристик $(l,w)$ город был разбит на районы сеткой переменного размера и в каждом районе выбрано по одной вершине-источнику и вершине-стоку корреспонденций. Размеры ячеек выбирались так, чтобы в каждом районе число людей, въезжающих и выезжающих из него было достаточно небольшим.
На рис.~\ref{fig:network} изображён дорожный граф центральной части Владивостока и его разбиение на районы.
Код и результаты моделирования размещены в репозитории \cite{kоtlyarova2020code}.
}

\color{black}

Авторы выражают благодарность проф. Ю.Е. Нестерову к 65-и летию которого приурочена данная статья, а также Мерузе Кубентаевой за постоянную помощь (консультации).

\dv{Авторы также выражают благодарность проф. Е.А. Нурминскогому за предоставленные по г. Владивостоку данные.}

Работа Е.В. Котляровой была выполнена в Сириусе (Сочи) в августе 2020 г. в рамках проектной студенческой смены <<Оптимизация, управление и информация>>.

Работа Е.В. Гасниковой была выполнена при поддержке Министерства науки и высшего образования Российской Федерации (госзадание) No. 075-00337-20-03, номер проекта 0714-2020-0005. Работа А.В. Гасникова была поддержана грантом РФФИ 18-29-03071 мк.


\end{document}